\documentclass{article}
\usepackage{latexsym}
\usepackage{amsmath,amsthm}
\usepackage{enumerate}
\usepackage{amssymb}
\usepackage{upgreek}
\usepackage{MnSymbol, wasysym}
\usepackage{graphics}
\usepackage{marvosym}
\usepackage{color}
\usepackage{color}

\def\f12{\frac 1 2}

\def\a{\alpha}

\def\pa{\partial}

\def\f12{\frac 1 2}

\newcommand{\lap}{\mbox{$\Delta \mkern-13mu /$\,}}

\newcommand{\nabb}{\mbox{$\nabla \mkern-13mu /$\,}}
\newtheorem{remark}{Remark}
\begin{document}

\title{A new physical-space approach to decay for the wave equation
with applications to black hole spacetimes}

\author{Mihalis Dafermos\thanks{Department of Pure Mathematics and Mathematical Statistics, University of Cambridge, Cambridge, CB3 0WB, United Kingdom}
 \and Igor Rodnianski\thanks{Department of Mathematics, Princeton University,
Fine Hall, Washington Road, Princeton, New Jersey 08544, United States}}

\maketitle
\begin{abstract}
We present a new general method for proving global decay of energy
through a suitable spacetime foliation, as well as
pointwise decay, starting from 
an integrated local energy decay estimate. The method is quite robust,
requiring only physical space techniques, and circumvents use of
multipliers or commutators with weights growing in $t$. 
In particular, the method applies to a wide class of perturbations
of Minkowski space as well as to Schwarzschild and Kerr black hole exteriors. 
\end{abstract}

\section{Introduction}
The wave equation
\begin{equation}
\label{THEWAVEEQ}
\Box_g\phi =0
\end{equation}
on general Lorentzian background metrics $g$ appears ubiquitously in mathematical physics,
often in the context of the linearisation of a  field
theory governed by a system of non-linear hyperbolic
p.d.e.'s. The most classical example occurs perhaps in the context of the Euler equations
of fluid mechanics, in which case $g$ represents the so-called acoustical metric. 
Another important source for $(\ref{THEWAVEEQ})$
arises from problems in general relativity, where $g$ represents the metric of spacetime,
and $(\ref{THEWAVEEQ})$ can be viewed as a poor-man's linearisation
of the Einstein equations.   For such metrics $g$, the global causal 
structure is often
much richer than that of  Minkowski space. In this context,
perhaps the most interesting examples 
are so-called ``black hole'' metrics. Briefly, these 
are metrics which--like Minkowski space--possess a natural
asymptotic structure representing future null infinity $\mathcal{I}^+$, but where the 
past of $\mathcal{I}^+$ has a non-trivial complement in the spacetime.

A fundamental problem for $(\ref{THEWAVEEQ})$ on suitable metrics
$g$ is to obtain ``decay estimates'', that is to 
say, estimates giving upper bounds for the rates of decay for $\phi$ toward the future,
in terms of quantities computable initially on a Cauchy surface. This
classical problem
can be studied as an end in itself, but it is often an essential step in proving stability results for
non-linear equations which give rise to $(\ref{THEWAVEEQ})$ as described above. From the latter
point of view, there is an important trade off between how much decay one proves and
how robust is the argument to prove it. The history of non-linear wave equations in the past
30 years, culminating in part in the monumental stability proofs (see for instance~\cite{book}) 
 based on robust estimates for
equations of the form $(\ref{THEWAVEEQ})$, are a testament to how important
it is to get this balance right.

The now standard approach to obtaining robust decay estimates proceeds by using weighted
multiplier estimates (going back to Morawetz~\cite{mora}) and weighted commutators
(going back to Klainerman~\cite{muchT}), both originating
in one form or another from the symmetries of Minkowski space.
In Minkowski space, full use of this method (applying both weighted multipliers
and weighted commutators) obtains (a) weighted energy decay, which can be expressed
for instance  as
\begin{eqnarray}
\label{fromZ}
\nonumber
\int_{t= \tau}&\big( (t^2/r^2+1) \phi^2+ (t+r+1)^2 ((\partial_{t}+\partial_{r})\phi)^2 \\
\nonumber
&+(|t-r|+1)^2((\partial_{t}-\partial_{r})\phi)^2 +  (t^2+r^2)|\nabb\phi|^2) \big)\\
 &\le C
\int_{t=0} (\phi^2+ r^2((\partial_t\phi)^2+(\partial_r\phi)^2+|\nabb\phi|^2),
\end{eqnarray}
where $\nabb$ denotes angular derivatives,
(b) pointwise decay for $\phi$ of the form
\begin{equation}
\label{fullMink}
|\phi|\le C \sqrt{E}(|t-r|+1)^{-1/2}(|t+r|+1)^{-1}
\end{equation}
for $t\ge 0$, where $E$ represents a higher order weighted energy at $t=0$, and
(c) pointwise estimates for derivatives of $\phi$, with additional decay, for instance:
\begin{equation}
\label{fullMinkder}
|\partial_t\phi|\le C \sqrt{E}(|t-r|+1)^{-3/2}(|t+r|+1)^{-1}.
\end{equation}
Obtaining (a) requires only use of inverted time translations
$(t^2+r^2)\partial_t+2tr\partial_r$ as a multiplier. Obtaining the full 
(b) and (c) requires both use of inverted time translations 
and commutation with the full algebra
of commutators associated to $\Box$.

Various applications to non-linear problems and other geometric settings 
have required variants of the vector field method where the use
of certain vector field multipliers and/or commutators
is avoided and less information is obtained concerning
(a), (b), or (c) above.  Klainerman 
and Sideris~\cite{KS} obtained decay results by commuting only with $t\partial_t+r\partial_r$
and angular momentum operators, and this, used together with integrated local decay estimates
in the style of Morawetz~\cite{mora2}, has proven useful 
in a number of problems concerning obstacles, multiple speeds, etc.~\cite{alinhac, age, ms}. A
common feature of all these approaches, however, is
the necessity of always applying at least some multipliers or 
commutators with weights growing in $t$.

In the present paper, we give a related but different approach to the problem of proving suitable
decay estimates for solutions to $(\ref{THEWAVEEQ})$ which, as we shall see,
has in principle a wider range of validity and may be even more 
useful for non-linear problems. In particular, the method here avoids both multipliers
and commutators with weights in $t$. 
Our approach does not start from scratch, but
begins with certain ``more basic'' estimates as ingredients, estimates which
have already been proven (either classically or much more recently)
in many cases of interest, by various methods. 
The method by which these ingredient estimates are obtained is in fact of no 
importance for the considerations here.
As we shall see, our approach applies to the study of $(\ref{THEWAVEEQ})$
on Minkowski space and
sufficiently small perturbations thereof, 
but also to Schwarzschild and Kerr black hole
spacetimes, and immediately yields essentially the decay consistent with
non-linear stability.
The common feature which these regimes share are 
\begin{enumerate}[(i)]
\item An integrated local energy decay estimate
\label{ide}
\item Good asymptotics towards null infinity allowing for a hierarchy of $r$-weighted energy identities
\item Uniform boundedness of energy
\label{iide}
\end{enumerate}
Estimates of the form~(\ref{ide}), already referred to above, were first
proven by Morawetz~\cite{mora2} for perturbations of Minkowski space,
while for Schwarzschild, estimates of the form~(\ref{ide}) are due to several authors
(see Section~4.4.1 of~\cite{jnotes}). 
For slowly rotating Kerr spacetimes $(\mathcal{M},g_{M,a})$ with parameters 
$|a|\ll M$, such estimates have been shown by~\cite{jnotes, tattoh} (see also~\cite{ab}), 
whereas,
for the general sub-extremal Kerr case $|a|<M$,
this will appear in a forthcoming paper. In the black hole case, the ``integrated decay''
estimate (\ref{ide}), to be useful here,
 must be understood in the sense of an estimate which does not degenerate 
at the horizon $\mathcal{H}^+$, 
an estimate which is indeed possible in view of a multiplier
construction intimately connected to the celebrated
red-shift effect (see~\cite{dr3, jnotes}). 
In the Schwarzschild and Kerr case, the integrated decay
estimate necessarily degenerates, however, near the analogue of the photon sphere,
and this loss is related to the trapping phenomenon; a non-degenerate
estimate can nevertheless  be obtained after losing differentiability. Similar phenomena
occur in  the simpler setting of obstacle problems in Minkowski space,
to which this method also applies.

We note that the hierarchy of $r$-weighted energy identities (ii) applied here corresponds to 
the use of a multiplier of the form $r^p(\partial_t+\partial_r)$ for $0\le p \le2$.

In principle, our approach can retrieve the ``full decay statements'' of the traditional
vector field method as would follow from the combined use
of~\cite{mora} and~\cite{muchT}, except that the energy decay statements
will in general lose derivatives if the integrated local energy decay statement does so.
In this short note,
we shall here only give the first step of the approach. This begins by 
obtaining $\tau^{-2}$ decay of the energy flux 
through a suitable foliation $\Sigma_\tau$ (whose leaves approach null infinity), 
where $\tau$ is the natural parameter of the foliation. 
This is precisely the analogue of~\eqref{fromZ}. (See the discussion of Section~\ref{refine}.)
This then leads immediately to pointwise results (in Section~\ref{pwise}) 
which are the analogue of
\[
|r^{1/2}\phi|\le Ct^{-1},\qquad |r\phi|\le Ct^{-1/2}.
\]
Again, no use is made of commutation
with vector fields with weights in $t$.
We outline in Section~\ref{refine2} 
how to  extend to the method to achieve the analogue of~\eqref{fullMink},~\eqref{fullMinkder}.

We consider the case of Minkowski space
first in Section~\ref{minkcase}, followed by Schwarzschild in Section~\ref{Scase}, 
to illustrate the structure of the argument
in each of these settings.  As shall be apparent, the nature of
the argument is such that elements (\ref{ide})--(\ref{iide}), suitably
interpreted, are sufficient and a very general theorem can in fact
be formulated. This general theorem includes in particular the Kerr case. 
See Section~\ref{genform}.
We leave the precise formulation of the most general assumptions for a follow-up paper.

\section{Notation}
Let $(\mathcal{M},g)$ be a Lorentzian manifold and $\phi$
a solution to $(\ref{THEWAVEEQ})$ on $\mathcal{M}$. 
We recall the energy-momentum tensor
\[
{\mathbb T}_{\mu\nu}[\phi]=\phi_{,\mu}\phi_{,\nu}-\frac12 g_{\mu\nu}\phi^{,\alpha}\phi_{,\alpha}.
\]
Given a vector field $X$, we define the currents
\[
J^X_\mu[\phi]= {\mathbb T}_{\mu\nu}[\phi]X^\nu, \qquad
K^X[\phi]= {\mathbb T}^{\mu\nu}[\phi]\pi^X_{\mu\nu},
\]
where $\pi^X_{\mu\nu}=\frac12 \mathcal{L}_Xg_{\mu\nu}$.
Recall that  
\[
\nabla^\mu J^X_\mu[\phi] = K^X[\phi].
\]
If $X$ is Killing, then $K^X[\phi]=0$,
in which case $J^X_\mu[\phi]$ is a conserved current.

Integrands without an explicit measure of integration 
are to be understood with respect to the induced volume form. 
For the case of null hypersurfaces $\Sigma$, integrands of the form 
$\int_\Sigma J^X_\mu n^\mu_\Sigma$ are to be understood with respect
to a choice of volume form on $\Sigma$ and corresponding normalisation for $n^\mu$
such that the integral represents the correct term in the energy identity. 
Alternatively, one can interpret this simply as the integral of the $3$-form
$\int_\Sigma {}^*J_X$,
where ${}^*$ denotes the Hodge star operator associated to the spacetime metric.

\section{The Minkowski case}
\label{minkcase}
Let us start by considering the wave equation in Minkowski space 
$$
\Box \phi=0.
$$ 
For convenience, we may assume for this discussion that $\phi$ is smooth and
compactly supported on a Cauchy hypersurface. This will ensure that all 
integrals are a priori defined. By a density argument, our results will only
require the finiteness of the quantity on the right hand side of the final estimate $(\ref{finalest})$
and the relation
$\phi\to 0$ at spacelike infinity.

Let $(t, r, \theta, \varphi)$ denote a choice of standard spherical polar coordinates 
on Minkowski space.  We shall let $\nabb$ denote the induced covariant derivative
on the spheres of constant $r$, and $\lap$ the induced Laplacian.
We will also define the null coordinates $u=t-r$, $v=t+r$.
Finally, let us introduce the notation $T$ to denote the coordinate vector field
$\partial_t$ with respect to $(t, r)$ coordinates.

For our argument it will be useful to consider the following foliation:
Fix $R>0$ and consider the hypersurfaces $\Sigma_\tau$ defined as the union of
the spacelike $\{t=\tau\}\cap\{r\le R\}$ and the null $\{u=(\tau-R)\}\cap \{ v\ge(\tau+R)\}$.

The following estimate can be shown by the original method of Morawetz~\cite{mora2}:
\begin{equation}
\label{Morawetz}
\tag{ILED-Mink}
\int_{\tau}^\infty \int_{r\le R} \left (J^T_\a[\phi] n^\a_{\Sigma_\tau} +\phi^2\right)
\le C_R \int_{\Sigma_\tau} J^T_\a[\phi] n^\a_{\Sigma_\tau}.
\end{equation}
The initials of the label stand for ``integrated local energy decay''.
Since $n^\alpha=T^\alpha$ in $r\le R$, we could replace the first term on the right
hand side more explicitly with the quantity $(\partial_t\phi)^2+(\partial_r\phi)^2+|\nabb\phi|^2$,
but we prefer this more geometric formulation to compare with other situations.
Our goal is to use \eqref{Morawetz} to derive decay for the energy flux
through $\Sigma_\tau$, and finally, pointwise decay estimates.

We rewrite the wave equation in the form 
$$
-\pa_u \pa_v \psi+\lap \psi=0,\quad \psi:=r\phi,
$$
multiply the equation by $r^p \pa_v\psi$ and integrate by parts in the region $D_{\tau_1}^{\tau_2}$ bounded by the
two null hypersurfaces $u_1=(\tau_1-R)$, $u_2=(\tau_2-R)$ and the timelike hypersurface 
$r=R$:
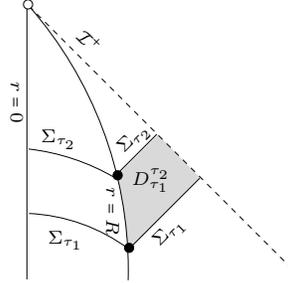
\[
\begin{picture}(0,0)%
\includegraphics{formink.pstex}%
\end{picture}%
\setlength{\unitlength}{2368sp}%
\begingroup\makeatletter\ifx\SetFigFont\undefined%
\gdef\SetFigFont#1#2#3#4#5{%
  \reset@font\fontsize{#1}{#2pt}%
  \fontfamily{#3}\fontseries{#4}\fontshape{#5}%
  \selectfont}%
\fi\endgroup%
\begin{picture}(3017,2953)(3371,-4819)
\put(4681,-3811){\makebox(0,0)[lb]{\smash{{\SetFigFont{7}{8.4}{\rmdefault}{\mddefault}{\updefault}{\color[rgb]{0,0,0}$D^{\tau_2}_{\tau_1}$}%
}}}}
\put(4966,-4462){\rotatebox{45.0}{\makebox(0,0)[lb]{\smash{{\SetFigFont{7}{8.4}{\rmdefault}{\mddefault}{\updefault}{\color[rgb]{0,0,0}$\Sigma_{\tau_1}$}%
}}}}}
\put(3801,-4421){\makebox(0,0)[lb]{\smash{{\SetFigFont{7}{8.4}{\rmdefault}{\mddefault}{\updefault}{\color[rgb]{0,0,0}$\Sigma_{\tau_1}$}%
}}}}
\put(3731,-3361){\makebox(0,0)[lb]{\smash{{\SetFigFont{7}{8.4}{\rmdefault}{\mddefault}{\updefault}{\color[rgb]{0,0,0}$\Sigma_{\tau_2}$}%
}}}}
\put(4596,-3482){\rotatebox{45.0}{\makebox(0,0)[lb]{\smash{{\SetFigFont{7}{8.4}{\rmdefault}{\mddefault}{\updefault}{\color[rgb]{0,0,0}$\Sigma_{\tau_2}$}%
}}}}}
\put(4090,-2266){\rotatebox{315.0}{\makebox(0,0)[lb]{\smash{{\SetFigFont{7}{8.4}{\rmdefault}{\mddefault}{\updefault}{\color[rgb]{0,0,0}$\mathcal{I}^+$}%
}}}}}
\put(3412,-2727){\rotatebox{270.0}{\makebox(0,0)[lb]{\smash{{\SetFigFont{7}{8.4}{\rmdefault}{\mddefault}{\updefault}{\color[rgb]{0,0,0}$r=0$}%
}}}}}
\put(4412,-3837){\rotatebox{270.0}{\makebox(0,0)[lb]{\smash{{\SetFigFont{7}{8.4}{\rmdefault}{\mddefault}{\updefault}{\color[rgb]{0,0,0}$r=R$}%
}}}}}
\end{picture}%

\]

We obtain
\begin{align}
\notag
&\int\limits_{u=\tau_2-R, v\ge \tau_2+R} r^p (\pa_v\psi)^2 \sin\theta\,d\theta\,
 d\phi\, dv \,\\
 \notag& +\int\limits_{D_{\tau_1}^{\tau_2}} r^{p-1}  \left (p(\pa_v\psi)^2 +
(2-p) |\nabb\psi|^2\right) \sin\theta\,d\theta\,
 d\phi\, du\, dv 
+\int\limits_{{\mathcal I}_{\tau_1-R}^{\tau_2-R}} r^p |\nabb\psi|^2
 \sin\theta\,d\theta\,
 d\phi \, du \notag\\ 
=& \int\limits_{u=\tau_1-R, v\ge \tau_1+R}r^p (\pa_v\psi)^2 \sin\theta\,d\theta\,
 d\phi\, dv  \notag \\
  &
+\int\limits_{\tau_1}^{\tau_2}  r^p \left (|\nabb\psi|^2- (\pa_v\psi)^2\right)\sin\theta\,d\theta\,
 d\phi\,  d\tau |_{r=R}\,.
\label{eq:en}\tag{p-WE-Mink}
\end{align}
The initials ``WE'' denote ``weighted energy''.
We have written explicitly the measure of integration to emphasize that the
expected $r^2$ part of the volume
form is now in fact included in the quadratic terms in $\psi$ (remember: $\psi=r\phi$). 
Thus, the factors $r^p$, $r^{p-1}$
are to be viewed as weights.
We see that the left hand side of the above identity is positive definite for $p\le 2$.

We note that identity~\eqref{eq:en} can be reinterpreted in terms of the energy identity
in $D^{\tau_2}_{\tau_1}$ satisfied by the currents
$J^V_\mu[\phi]$, $K^V[\phi]$ corresponding to $V=r^p\partial_v$ (where $\partial_v$ is understood
with respect to $(u,v)$ coordinates), suitably modified, however,
 by appropriate $0$'th order terms.

Let us first apply \eqref{eq:en} with $p=2$. 
Then, observing that the last term on the right hand side of \eqref{eq:en}
can be controlled (after a bit of averaging in $R$) by the left hand side of \eqref{Morawetz}, 
we obtain
\begin{eqnarray}
\label{afterav}
\nonumber
\int\limits_{u=\tau_2-R, v\ge \tau_2+R} r^2 (\pa_v\psi)^2
\sin\theta\,d\theta\,
 d\phi\,  dv +\int\limits_{D_{\tau_1}^{\tau_2}} r 
(\pa_v\psi)^2 \sin\theta\,d\theta\,
 d\phi\,  du\, dv \\
 \le \int\limits_{u=\tau_1-R, v\ge \tau_1+R}r^2 (\pa_v\psi)^2
 \sin\theta\,d\theta\,d\phi\,  dv +  
C\int_{\Sigma_{\tau_1}} J^T_\a[\phi] n^\a_{\Sigma_\tau}.
\end{eqnarray}
This implies that we can find a dyadic sequence of $\tau_n\to \infty$ with the property that 
\begin{eqnarray*}
\int\limits_{u=\tau_n-R, v\ge \tau_n+R} r (\pa_v\psi)^2
\sin\theta\,d\theta\, d\phi\,  dv
\\
\le C \tau_n^{-1} \left [\,\,\int\limits_{u=\tau_1-R, v\ge \tau_1+R}r^2 (\pa_v\psi)^2
\sin\theta\,d\theta\,
 d\phi \, dv +  
\int_{\Sigma_{\tau_1}} J^T_\a[\phi] n^\a_{\Sigma_\tau}\right].
\end{eqnarray*}

We now apply \eqref{eq:en} with $p=1$ to the region $D_{\tau_{n-1}}^{\tau_n}$ to obtain 
\begin{align*}
&\int\limits_{u=\tau_n-R, v\ge \tau_n+R} r (\pa_v\psi)^2 
\sin\theta\,d\theta\,
 d\phi\, dv +\int\limits_{D_{\tau_{n-1}}^{\tau_{n}}} 
\left ((\pa_v\psi)^2 +
 |\nabb\psi|^2\right) \sin\theta\,d\theta\,
 d\phi\,  du \,dv
\\ &\hskip5pc\le 
\int\limits_{u=\tau_{n-1}-R, v\ge \tau_{n}+R} r (\pa_v\psi)^2 
\sin\theta\,d\theta\,
 d\phi\, dv +
C\int_{\Sigma_{\tau_{n-1}}} J^T_\a[\phi] n^\a_{\Sigma_\tau}\\ &\hskip 5pc \le
C \tau_n^{-1} \left [\,\,\int\limits_{u=\tau_1-R, v\ge \tau_1+R}r^2 (\pa_v\psi)^2 
\sin\theta\,d\theta\,
 d\phi\,  dv +  
\int_{\Sigma_{\tau_1}} J^T_\a[\phi] n^\a_{\Sigma_\tau}\right]\\ &\hbox{}\hskip5pc+  C
\int_{\Sigma_{\tau_{n-1}}} J^T_\a[\phi] n^\a_{\Sigma_\tau}.
\end{align*}
Adding a multiple of the estimate \eqref{Morawetz}, we obtain
\begin{align}
&\int_{\tau_{n-1}}^{\tau_n} \int_{r\le R} \left (J^T_\a[\phi] n^\a_{\Sigma_\tau} +\phi^2\right)
\notag
\\
& +
\int\limits_{u=\tau_n-R, v\ge \tau_n+R} r (\pa_v\psi)^2 
\sin\theta\,d\theta\,
 d\phi\, dv +\int\limits_{D_{\tau_{n-1}}^{\tau_{n}}} 
\left ((\pa_v\psi)^2 +
 |\nabb\psi|^2\right) \sin\theta\,d\theta\,
 d\phi\,  du\, dv
\notag\\ &\hskip5pc \le
C \tau_n^{-1} \left [\,\,\int\limits_{u=\tau_1-R, v\ge \tau_1+R}r^2 (\pa_v\psi)^2
\sin\theta\,d\theta\,
 d\phi\,  dv +  
\int_{\Sigma_{\tau_1}} J^T_\a[\phi] n^\a_{\Sigma_\tau}\right]
\notag
\\ &\hbox{}\hskip 5pc+  C
\int_{\Sigma_{\tau_{n-1}}} J^T_\a[\phi] n^\a_{\Sigma_\tau}.
\label{eq:1}
\end{align}
We now observe that 
$$
\int_{v_n}^\infty \left(\pa_v (r\phi)\right)^2 dv = \int_{v_n}^\infty \left [r^2 (\pa_v \phi)^2 + 2r \pa_v\phi\,\phi +\phi^2\right]dv = \int_{v_n}^\infty r^2 (\pa_v \phi)^2 dv  - r \phi^2|_{v=v_n}.
$$
Substituting this into \eqref{eq:1} we obtain
\begin{align}
\label{biv}
\nonumber
&\int_{\tau_{n-1}}^{\tau_n} \int_{\Sigma_\tau} J^T_\a[\phi] n^\a_{\Sigma_\tau}  \le
C \tau_n^{-1} \left [\,\,\int\limits_{u=\tau_1-R, v\ge \tau_1+R}r^2 (\pa_v\psi)^2 
\sin\theta\,d\theta\,
 d\phi\,  dv +  
\int_{\Sigma_{\tau_1}} J^T_\a[\phi] n^\a_{\Sigma_\tau}\right]\\ &\hskip 5pc+  C
\int_{\Sigma_{\tau_{n-1}}} J^T_\a[\phi] n^\a_{\Sigma_\tau}.
\end{align}

Finally, in view  also of  the energy bound
\begin{equation}
\label{eb}
\tag{EB-Mink}
\int_{\Sigma_\tau}J^T_\a[\phi] n^\a_{\Sigma_\tau}\le
\int_{\Sigma_{\tau'}} J^T_\a[\phi] n^\a_{\Sigma_{\tau'}},
\end{equation}
which holds for all $\tau\ge \tau'$,
the estimate $(\ref{biv})$ easily implies that 
\begin{equation}
\label{finalest}
\int_{\Sigma_{\tau}} J^T_\a[\phi] n^\a_{\Sigma_\tau}\le
C \tau^{-2} \left [\,\,\int\limits_{u=\tau_1-R, v\ge \tau_1+R}r^2 (\pa_v\psi)^2
\sin\theta\,d\theta\,
 d\phi\,  dv +  
\int_{\Sigma_{\tau_1}} J^T_\a[\phi] n^\a_{\Sigma_\tau}\right]
\end{equation}
for any $\tau\ge \tau_1$.
\vskip1pc
\begin{remark}
Note that $(\ref{finalest})$
essentially corresponds to $(\ref{fromZ})$ of the discussion of the vector field method.
In view of the absense of weights in $t$ in the multiplier used for the proof, it
is already immediate from the above that $(\ref{finalest})$ follows 
similarly on sufficiently small non-stationary perturbations of Minkowski space. 
This estimate is in itself very useful for applications to
global existence problems for quasilinear equations, and previously
was not available without commutation with inverted time translations with 
their quadratic weights in $t$.
See the discussion in Section~\ref{refine}.
\end{remark}

\section{The Schwarzschild case}
\label{Scase}
We now consider the wave equation $(\ref{THEWAVEEQ})$
on a Schwarzschild background $(\mathcal{M},g)$ with parameter $M$.
See Section~2 of~\cite{jnotes}.
Let $\mathcal{D}$ denote the closure of (a connected component of)
the domain of outer communications, 
let $\Sigma_{\tau}$ be a translation-invariant family of hypersurfaces
connecting the event horizon $\mathcal{H}^+$ and null infinity $\mathcal{I}^+$ 
as below:
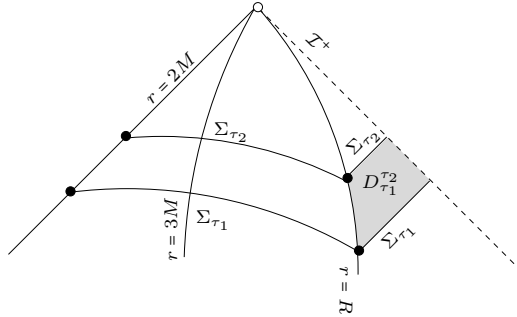
\[
\begin{picture}(0,0)%
\includegraphics{forschw.pstex}%
\end{picture}%
\setlength{\unitlength}{2368sp}%
\begingroup\makeatletter\ifx\SetFigFont\undefined%
\gdef\SetFigFont#1#2#3#4#5{%
  \reset@font\fontsize{#1}{#2pt}%
  \fontfamily{#3}\fontseries{#4}\fontshape{#5}%
  \selectfont}%
\fi\endgroup%
\begin{picture}(5419,3237)(969,-5103)
\put(4681,-3811){\makebox(0,0)[lb]{\smash{{\SetFigFont{7}{8.4}{\rmdefault}{\mddefault}{\updefault}{\color[rgb]{0,0,0}$D^{\tau_2}_{\tau_1}$}%
}}}}
\put(4966,-4462){\rotatebox{45.0}{\makebox(0,0)[lb]{\smash{{\SetFigFont{7}{8.4}{\rmdefault}{\mddefault}{\updefault}{\color[rgb]{0,0,0}$\Sigma_{\tau_1}$}%
}}}}}
\put(4596,-3482){\rotatebox{45.0}{\makebox(0,0)[lb]{\smash{{\SetFigFont{7}{8.4}{\rmdefault}{\mddefault}{\updefault}{\color[rgb]{0,0,0}$\Sigma_{\tau_2}$}%
}}}}}
\put(4452,-4627){\rotatebox{270.0}{\makebox(0,0)[lb]{\smash{{\SetFigFont{7}{8.4}{\rmdefault}{\mddefault}{\updefault}{\color[rgb]{0,0,0}$r=R$}%
}}}}}
\put(4090,-2266){\rotatebox{315.0}{\makebox(0,0)[lb]{\smash{{\SetFigFont{7}{8.4}{\rmdefault}{\mddefault}{\updefault}{\color[rgb]{0,0,0}$\mathcal{I}^+$}%
}}}}}
\put(2483,-2946){\rotatebox{45.0}{\makebox(0,0)[lb]{\smash{{\SetFigFont{7}{8.4}{\rmdefault}{\mddefault}{\updefault}{\color[rgb]{0,0,0}$r=2M$}%
}}}}}
\put(3131,-3261){\makebox(0,0)[lb]{\smash{{\SetFigFont{7}{8.4}{\rmdefault}{\mddefault}{\updefault}{\color[rgb]{0,0,0}$\Sigma_{\tau_2}$}%
}}}}
\put(2941,-4181){\makebox(0,0)[lb]{\smash{{\SetFigFont{7}{8.4}{\rmdefault}{\mddefault}{\updefault}{\color[rgb]{0,0,0}$\Sigma_{\tau_1}$}%
}}}}
\put(2749,-4580){\rotatebox{90.0}{\makebox(0,0)[lb]{\smash{{\SetFigFont{7}{8.4}{\rmdefault}{\mddefault}{\updefault}{\color[rgb]{0,0,0}$r=3M$}%
}}}}}
\end{picture}%

\]
For definiteness, one can choose the family $\Sigma_\tau$ explicitly
as follows: Given a constant $R> 3M$,
first define a Regge-Wheeler coordinate $r^*= r+2M\log (r-2M)-R$,
the null coordinates $u= t-r^*$, $v=t+r^*$ and the coordinate
$t^*= t+2M\log(r-2M)$. Then $\Sigma_\tau$
can be defined to coincide with $t^*=\tau$ for $r\le R$, and to coincide with $u=\tau$
for $r\ge R$. Let us define $D^{\tau_2}_{\tau_1}$ analogously as before,
i.e., with the present normalisation of the null coordinates: 
$D^{\tau_2}_{\tau_1}=\{r\ge R\} \cap\{\tau_1\le u\le \tau_2\}$.

It will be useful to fix a translation-invariant timelike vector field $N$ on $J^+(\Sigma_{\tau_1})$
such that $N=\partial_t $ for $r\ge R$, say. (The coordinate vector field $\partial_{t^*}$ in
$(r, t^*)$ coordinates will do for instance.) 
Let us note that~\eqref{eb} is replaced by
\begin{equation}
\label{eb-S}
\tag{EB-Schw}
\int_{\Sigma_\tau}J^N_\a[\phi] n^\a_{\Sigma_\tau}\le
C\int_{\Sigma_{\tau'}} J^N_\a[\phi] n^\a_{\Sigma_{\tau'}}.
\end{equation}
This non-degenerate energy boundedness statement 
was originally proven in~\cite{dr3}
but can in fact be shown independently of decay results, as in~\cite{bounded} 
and Section~3 of~\cite{jnotes}.

The analogue of~\eqref{Morawetz} is 
\begin{equation}
\label{M1Sc}
\tag{ILED1-Schw}
\int_{\tau}^\infty \int_{r\le R} \left (\chi J^N_\a[\phi] n^\a_{\Sigma_\tau} +(\partial_r\phi)^2 +\phi^2\right)
\le C \int_{\Sigma_\tau} J^N_\a[\phi] n^\a_{\Sigma_\tau},
\end{equation}
where $\chi$ is a weight which vanishes quadratically at $r=3M$.
This type of estimate, which does not degenerate on the horizon,
was originally proven in~\cite{dr3}. For background, see Section 4.4.1 of~\cite{jnotes}.
Commuting with $T$, we may obtain for instance
\begin{equation}
\label{M2Sc}
\tag{ILED2-Schw}
\int_{\tau}^\infty \int_{r\le R} J^N_\a[\phi] n^\a_{\Sigma_\tau} 
\le C \int_{\Sigma_\tau} J^N_\a[T \phi] n^\a_{\Sigma_\tau}+
 C\int_{\Sigma_\tau} J^N_\a[\phi] n^\a_{\Sigma_\tau}.
\end{equation}

We will outline more explicitly further on in this section
how the above are used in the Schwarzschild
case. Let us immediately remark, however, that~\eqref{M1Sc} would
be sufficient for the use of \eqref{Morawetz} to obtain the analogue of
$(\ref{afterav})$, whereas one would have to use~\eqref{M2Sc} to obtain~\eqref{eq:1},
necessitating the appearance of a higher order quantity on the right hand side.
This will mean that, in the Schwarzschild case,
the final estimate loses differentiability, as 
expected.\footnote{The 
precise loss of differentiability in~\eqref{M2Sc} is wasteful.
Using a refined estimate (see~\cite{marzuola})
would lead to an improved result here in the sense of loss of differentiability, but
this would take us outside the realm of physical space methods.} In this latter use of~\eqref{M2Sc},
it is essential that one has the non-degenerate quantity $J^N[\phi]$ on the horizon.

Before adapting the argument of Section~\ref{minkcase} to Schwarzschild, it remains to 
derive the analogue of \eqref{eq:en}.

As before, we write the wave equation $(\ref{THEWAVEEQ})$ as 
$$
-\pa_u\pa_v \psi + (1-\frac {2M}r)\lap \psi-\frac{2M(1-\frac{2M}r)}{r^3} \psi=0,
\qquad \psi :=r\phi.
$$
We may rewrite this in the form 
$$
-\frac 1{1-\frac{2M}r}\pa_u\pa_v \psi + \lap \psi-\frac{2M}{r^3} \psi=0.
$$
Observe that this implies that
\begin{eqnarray*}
-\frac 12 \pa_u \left(\frac {r^p}{1-\frac{2M}r} |\pa_v\psi|^2\right)+ \frac 12 \pa_u \left(\frac {r^p}{1-\mu}\right) |\pa_v\psi|^2\\
-
\frac 12 \pa_v\left (\frac{2M r^p}{r^3} \psi^2\right)+\frac 12 \pa_v\left (\frac{2M r^p}{r^3}\right) \psi^2 +
r^p\lap \psi\pa_v \psi=0,
\end{eqnarray*}
which generates the following additional terms in the analogue of~\eqref{eq:en}:
a boundary term at null infinity of the right sign, proportional to
$$
2M \int\limits_{{\mathcal I}_{\tau_1-R}^{\tau_2-R}} r^{p-3}  |\psi|^2 
\sin\theta\,d\theta\,
 d\phi\,du 
$$
and the bulk terms 
\begin{align*}
&-M \int\limits_{D_{\tau_1}^{\tau_2}} \pa_v \left [r^{p-3}\right] \psi^2  \sin\theta\,d\theta\,
 d\phi\,du \, dv =
(3-p) M \int\limits_{D_{\tau_1}^{\tau_2}} r^{p-4} \psi^2  \sin\theta\,d\theta\,
 d\phi\,du \, dv  \\&
-\frac 12\int\limits_{D_{\tau_1}^{\tau_2}} \pa_u \left [\frac {r^{p}}{1-\mu}\right] |\pa_v\psi|^2 
 \sin\theta\,d\theta\,
 d\phi\,du \, dv =
\int\limits_{D_{\tau_1}^{\tau_2}} r^{p-1}\left [p-
\frac {2M}{r}\right] |\pa_v\psi|^2 \sin\theta\,d\theta\,
 d\phi\,du \, dv  .
\end{align*}
Both expressions have the right sign as long as $0<p\le 3$ and $r$ is sufficiently large.
Note also the new $(1-2M/r)$ weights in the analogue of~\eqref{eq:en},
which of course will play no role for large $r$.

Suppressing the spherical integration from
the notation, we thus obtain, for $R$ sufficiently large, the estimate 
\begin{align}
\notag
&\int\limits_{u=\tau_2, v\ge \tau_2+R} r^p (\pa_v\psi)^2 dv +\int\limits_{D_{\tau_1}^{\tau_2}} r^{p-1}  \left (p(\pa_v\psi)^2 +
(2-p) |\nabb\psi|^2\right) du dv
+\int\limits_{{\mathcal I}_{\tau_1-R}^{\tau_2-R}} r^p |\nabb\psi|^2\notag\\ &\hskip 3pc
\le C\left(\int\limits_{u=\tau_1-R, v\ge \tau_1+R}r^p (\pa_v\psi)^2 dv
+\int\limits_{\tau_1}^{\tau_2}  r^p \left (|\nabb\psi|^2+ (\pa_v\psi)^2\right) d\tau |_{r=R}\right).
\label{eq:Sc}\tag{p-WE-Schw}
\end{align}

Given \eqref{M1Sc}, \eqref{M2Sc}, \eqref{eq:Sc}, \eqref{eb-S},
one repeats the Minkowski argument
until before $(\ref{eq:1})$.

At this point, one adds the estimate~\eqref{M2Sc} to obtain
\begin{align}
&\int_{\tau_{n-1}}^{\tau_n} \int_{r\le R} \left (J^N_\a[\phi] n^\a_{\Sigma_\tau} +\phi^2\right) +
\int\limits_{u=\tau_n-R, v\ge \tau_n+R} r (\pa_v\psi)^2 dv +\int\limits_{D_{\tau_{n-1}}^{\tau_{n}}} 
\left ((\pa_v\psi)^2 +
 |\nabb\psi|^2\right) du dv
\notag\\ &\hskip5pc \le
C \tau_n^{-1} \left [\,\,\int\limits_{u=\tau_1-R, v\ge \tau_1+R}r^2 (\pa_v\psi)^2 dv +  
\int_{\Sigma_{\tau_1}} J^N_\a[\phi] n^\a_{\Sigma_\tau}\right]\notag\\ &\hbox{}\hskip 5pc+  C
\int_{\Sigma_{\tau_{n-1}}} J^N_\a[\phi] n^\a_{\Sigma_\tau}
+  C
\int_{\Sigma_{\tau_{n-1}}} J^N_\a[T\phi] n^\a_{\Sigma_\tau}\label{eq:1Sc}.
\end{align}
One obtains as before
\begin{align*}
&\int_{\tau_{n-1}}^{\tau_n} \int_{\Sigma_\tau} J^N_\a[\phi] n^\a_{\Sigma_\tau}  \le
C \tau_n^{-1} \left [\,\,\int\limits_{u=\tau_1-R, v\ge \tau_1+R}r^2 (\pa_v\psi)^2 dv +  
\int_{\Sigma_{\tau_1}} J^N_\a[\phi] n^\a_{\Sigma_\tau}\right]\\ &\hskip 10pc+  C
\int_{\Sigma_{\tau_{n-1}}} J^N_\a[\phi] n^\a_{\Sigma_\tau}
+ C
\int_{\Sigma_{\tau_{n-1}}} J^N_\a[T\phi] n^\a_{\Sigma_\tau}
\\ & \hskip4pc \le
C\tau^{-1}_n\int_{\Sigma_{\tau_{1}}} J^N_\a[\phi] n^\a_{\Sigma_\tau}
 +C
\int_{\Sigma_{\tau_{1}}} J^N_\a[\phi] n^\a_{\Sigma_\tau}
+ C
\int_{\Sigma_{\tau_{1}}} J^N_\a[T\phi] n^\a_{\Sigma_\tau}.
\end{align*}

Identifying a good dyadic sequence, applying the boundedness statement~\eqref{eb-S}, 
and using the above inequality again with $\phi$ replaced by $T\phi$,
we obtain finally
\begin{equation}
\label{f-S}
\int_{\Sigma_{\tau}} J^N_\a[\phi] n^\a_{\Sigma_\tau}\le
C \tau^{-2}\sum_{i=0}^2
\int_{\Sigma_{\tau_1}} r^2J^N_\a[T^i\phi] n^\a_{\Sigma_\tau}.
\end{equation}

\section{Pointwise estimates}
\label{pwise}
Starting from the energy decay bounds $(\ref{finalest})$, $(\ref{f-S})$ proven above,
one can obtain pointwise decay estimates
after additional commutations and applications of weighted Sobolev inequalities.
To obtain decay near null infinity, we shall require commutator vector fields
with weights in $r$,
but, again in accordance with the philosophy of our method, no weights in $t$ will be used. 

In the Schwarzschild and Minkowski space, one approach for this is to use
only commutation with the so-called angular momentum operators $\Omega_i$, defined
as  a  standard basis for the Lie algebra corresponding to rotations. 
See~\cite{dr3}.

An alternative argument (see Section~5.3.8 of~\cite{jnotes}) 
uses commutation with cut-off versions of $\Omega_i$, supported only
near infinity, together with additional  commutations with $T$ and the red-shift
vector field $N$ (or equivalently $Y$) and an appeal to elliptic estimates. 
 This argument is more robust, and moreover,
can also be used to obtain non-degenerate pointwise estimates of higher derivatives up
to and including the horizon.

From this latter approach, we obtain finally:
\begin{equation}
\label{ourpw}
|r^{1/2} \phi| \le C\sqrt{E}\tau^{-1}, \qquad
|r\phi|\le C\sqrt{E}\tau^{-1/2},
\end{equation}
where 
\[
E= \sum_{|\alpha|\le 2}\sum_{\Gamma = T, N, \Omega_i}
\int_{\Sigma_0} r^2(J_\mu^N(\Gamma^{\alpha}\phi)+J_\mu^N(\Gamma^\alpha T\phi)
+J_\mu^N(\Gamma^\alpha TT\phi)n^\mu_{\Sigma_0}.
\]
We note  that the second two terms of the integrand arise from the loss
in \eqref{M2Sc} and are thus unnecessary in the
case of Minkowski space, where, in addition, of course, we take $N=T$.

\section{A general formulation and the Kerr case}
\label{genform}
It should be clear from the structure of the above argument that the properties 
(i), (ii), (iii), as captured by suitable estimates
 in the form of  \eqref{M1Sc}--\eqref{M2Sc}, 
\eqref{eq:Sc}, \eqref{eb-S}, are in fact sufficient for 
obtaining pointwise-in-time energy decay of the form  $(\ref{f-S})$. 
In particular, we note that the stationarity of the metric is not used directly,
as long as one has the above estimates (i), (ii), (iii),
and (in the case where there is trapping),
as long as one can absorb the errors arising from commuting with $T$. 
We shall give a general formulation in a follow-up paper. 
The results apply to a wide class of background metrics and can include
obstacles with trapped rays, black holes, etc.

Of particular importance are the Kerr spacetimes $(\mathcal{M},g_{a,M})$. See Section~5.1
of~\cite{jnotes} for a discussion. In the case $|a|\ll M$, the estimate associated to (ii) 
follows by the  results of~\cite{jnotes, tattoh}, whereas
the estimate associated to (iii) follows from~\cite{bounded}.
For the general case $|a|<M$, the estimates associated to (i) and (ii) will follow
by forthcoming results.
On the other hand, a computation reveals that the 
analogue of (i) follows as in the Schwarzschild case. 
Thus, one obtains precisely $(\ref{f-S})$ for Kerr, where $T$ is the stationary Killing field,
$N$ is a globally timelike translation-invariant vector field, and $\Omega_i$
are angular momentum operators on an ambient background Schwarzschild metric.

In the general formulation,
to obtain pointwise estimates as in $(\ref{ourpw})$,
the assumptions necessary for $(\ref{f-S})$ have to be supplemented with
assumptions regarding the existence of sufficiently many good commutators. 
Again, we leave the most general formulation to a follow-up paper. We remark here,
however, that the necessary properties can be immediately explicitly verified for Kerr,
for the entire range $|a|<M$, in view of the
timelike nature of the span of $\partial_t$ and $\partial_\phi$
outside the horizon, and the positivity of the surface gravity  ensuring the validity
of the commutation Theorem~7.2 of~\cite{jnotes}.
Thus, as in Schwarzschild, one obtains the pointwise~$(\ref{ourpw})$ for Kerr, where
the vector fields are interpreted as in the previous paragraph.

\section{Discussion}
\label{refine}
As mentioned before, the energy decay estimate $(\ref{finalest})$ 
corresponds exactly to what one can obtain in Minkowski space
by applying the multiplier $J^Z_\mu$ associated to the generator
of inverted time translations
\[
Z= u^2\partial_u+v^2\partial_v=(t^2+r^2)\partial_t+2tr\partial_r.
\]
This argument for obtaining decay estimates in Minkowski space
goes back to Morawetz~\cite{mora}.
Analogues of the current $J^{Z}_\mu$ 
can be constructed in the Schwarzschild case~\cite{dr3, bs} and
slowly rotating Kerr case $|a|\ll M$~\cite{jnotes, ab}, and were used
in conjunction with an integrated
decay type inequality, to obtain decay of energy in $\Sigma_\tau$. Already, however,
in the Kerr case, there is a loss in the power of decay which one can obtain
by this method, which is in principle proportional to the parameter $a$. 
See Section~5.3.6 of~\cite{jnotes}. The difficulty with these 
constructions arises from error terms in the region $r\le R$ which inherit large weights in $t$.
These problems could be even more profound in the setting of non-linear  applications.
This issue was in fact our motivation for the present work.

It is amusing to compare the role of null infinity $\mathcal{I}^+$ in the present argument
with the role of the cosmological horizon in the Schwarzschild-de Sitter case studied in~\cite{dr4}. 
There, the analogue of the iteration process applied here 
could be continued indefinitely to yield decay bounds
faster than any polynomial rate. The reason is that
in the analogue of \eqref{eq:Sc}, the weights of the boundary terms coincide with the weights of
the bulk terms. This is precisely related to the red-shift.  Indeed, this is reflected already in 
the Schwarzschild case 
at the event horizon by the fact that $J^N[\phi]$ appears without degeneration at $r=2M$
on both sides of inequality~\eqref{M1Sc}. A similar relation holds at the cosmological horizon of
Schwarzschild-de Sitter. In contrast, in the present asymptotically flat case, 
the weights in \eqref{eq:Sc} are polynomial in $r$ and drop by $1$ power when
comparing bulk and boundary terms.
This is why it is important to keep the whole $p$-hierarchy
of estimates \eqref{eq:Sc} in mind for the method to be useful.

We mention finally a forthcoming Fourier-based method of Tataru and collaborators for 
obtaining pointwise estimates
for stationary metrics starting again from (i) and good asymptotic behaviour of the 
metric\footnote{talk of D.~Tataru, Oberwolfach, October 2009}.

\section{Refinements: the full decay of the vector field method}
\label{refine2}
As explained in the introduction, the ``full vector field method'' yields faster pointwise
decay of the form $(\ref{fullMink})$.
One way to prove such a result starting
from our $(\ref{f-S})$ would be to apply one commutation with the scaling
vector field $t\partial_t+r\partial_r$. (See work of Luk~\cite{luk} for the Schwarzschild
case.)
This would be contrary, however, to 
the philosophy of our current approach, which is to avoid any multipliers and
commutators with weights in $t$.
It turns out that our approach can in fact be extended to yield more decay
without such commutations, by instead 
commuting with $\partial_v$ (again understood in a null $(u,v)$ coordinate system),
and $\Omega_i$
which allows applying the hierarchy of estimates \eqref{eq:en} or~\eqref{eq:Sc} 
for larger weights in $r$,
i.e.~for larger values of the parameter $p$. 
This leads
to additional
decay for a higher-order energy flux through $\Sigma_\tau$, which can then 
be transformed into the usual 
decay for $\phi$ using also ideas of~\cite{luk}. 
This will be discussed in a forthcoming follow-up paper.

\section*{Acknowledgments}
This work was conducted when I.R.~was visiting Cambridge in February 2009. 
I.R.~thanks the University of Cambridge for hospitality. M.D.~is supported in part by
a grant from the European Research Council. I.R.~is supported in part by NSF
grant DMS-0702270.

\end{document}